\documentclass{amsart}
\pdfoutput=1
\usepackage{caption}
\usepackage{subcaption}
\usepackage{amssymb}
\usepackage{graphicx}
\usepackage{wrapfig}
\begin{document}

\author{G\'abor Fejes T\'oth}
\address{Alfr\'ed R\'enyi Institute of Mathematics,
Re\'altanoda u. 13-15., H-1053, Budapest, Hungary}
\email{gfejes@renyi.hu}

\author{L\'aszl\'o Fejes T\'oth}

\author{W{\l}odzimierz. Kuperberg}
\address{Department of Mathematics \& Statistics, Auburn University, Auburn, AL36849-5310, USA}
\email{kuperwl@auburn.edu}

\title{Ball packings in hyperbolic space}
\thanks{The English translation of the book ``Lagerungen in der Ebene,
auf der Kugel und im Raum" by L\'aszl\'o Fejes T\'oth will be
published by Springer in the book series Grundlehren der
mathematischen Wissenschaften under the title
``Lagerungen---Arrangements in the Plane, on the Sphere and
in Space". Besides detailed notes to the original text the
English edition contains eight self-contained new chapters
surveying topics related to the subject of the book but not
contained in it. This is a preprint of one of the new chapters.}

\begin{abstract}
In hyperbolic space density cannot be defined by a limit as we
define it in Euclidean space. We describe the local density
bounds for sphere packings and we discuss the different attempts
to define optimal arrangements in hyperbolic space.
\end{abstract}

\maketitle

It is natural to extend the study of packing and covering problems to the
hyperbolic plane, as well as hyperbolic spaces of higher dimension.
However, we encounter here the problem of defining a ``reasonable''
notion of global density. Since in hyperbolic space the volume and
surface area of a ball of radius $r$ are of the tsame order of magnitude,
density cannot be defined by a limit as we define it in Euclidean space.
{\sc{B\"or\"oczky}} \cite{Boroczky74} exposed that the problem is much
deeper than one would expect. Namely, he constructed an arrangement of
congruent circles along with two different decompositions $Z_1$ and $Z_2$
of the hyperbolic plane into congruent cells, each containing one circle,
with the property that the circles' density in each cell of $Z_i$ presents
the same value $d_i$ ($i=1,2$), and yet $d_1\neq d_2$.

Two ways were considered to overcome the aforementioned difficulties with
defining global density in hyperbolic geometry. One way was to consider local
density relative to Dirichlet cells; the other was to come up with suitable notions
alternate to arrangements of optimal density.
\section{The simplex bound}
Recall the simplex bound of {\sc{B\"or\"oczky}} \cite{Boroczky78}:
For any packing of balls of radius $r$ in $n$-dimensional
hyperbolic space the density of a ball in its Dirichlet cell cannot exceed
the density $d_n(r)$ of $n+1$ mutually touching balls of radius $r$ with
respect to the simplex spanned by the centers of the balls. B\"or\"oczky
stated this for balls of finite radius, however the proof can be extended
to packings of horoballs (balls of infinite radius). As in the case of balls
of finite radius, the bound refers to the density of the horoballs in their
Dirichlet cell in this case as well. It should be mentioned that the volume
of the horoballs, as well as their Dirichlet cells, is infinite. However, as
the hyperbolic metric on a horosphere is Euclidean, the density of a horoball
in its Dirichlet cell can be defined by a limit of the density in growing
sectors of the horoball.

The definition of the simplex bound also needs further explanation in the case
of horoball packings. Let $S$ be a totally asymptotic regular simplex in
$n$-dimensional hyperbolic space. Take $n+1$ mutually tangent horoballs
centered at the vertices of $S$ and consider their density relative to $S$.
In contrast to balls of finite radius, the condition that the horoballs are mutually
tangent does not determine the configuration. We get the simplex bound $d_n(\infty)$
occurring in B\"or\"oczky's theorem if the symmetry group of the arrangement
of the horoballs coincides with that of $S$.

{\sc{Szirmai}} \cite{Szirmai12,Szirmai13} investigated the density of mutually
tangent horoballs centered at the vertices of a totally asymptotic regular
simplex $S$ in the case when the symmetry group of the arrangement of the
horoballs does not coincides with that of $S$. It turned out that for $n\ge3$
there are arrangement of horoballs that produce local density in $S$ higher
than $d_n(\infty)$. Of course, these configurations are only locally optimal
and cannot be extended to the entire hyperbolic space. But their existence
shows that there might also be arrangements of mutually tangent balls of finite
radius $r$ centered at the vertices of a simplex with density in the simplex
higher than $d_n(r)$.

We recall that for a packing of circles of radius $r$ in a plane of
constant curvature the density of each circle in its Dirichlet cell, as well
the density of the circles in each Delone cell, is at most $d_2(r)$. It
appears that in higher-dimensional hyperbolic space $d_n(r)$ is an upper
bound only for the local density in the Dirichlet cells and not for
the density in the Delone simplices, although no explicit example is
known. Also, it is an interesting open question whether, for some $n\ge3$,
in $n$-dimensional Euclidean or spherical space the density of a packing
of balls of radius $r$ in the Delone simplices can exceed $d_n(r)$.

In 3-dimensional hyperbolic space there is a remarkable tiling $\{6,3,3\}$
whose cells are degenerate Euclidean polyhedra $\{6,3\}$, circumscribed about
horoballs. The horoballs inscribed in the cells of this tiling form a packing
while the horoballs circumscribed about the cells form a covering.
{\sc{Coxeter}} \cite{Coxeter54} calculated the density of the horoballs
in the cells and found that it is
$$
d_3(\infty)=\left(
1+\frac{1}{2^2}-\frac{1}{4^2}-\frac{1}{5^2}+\frac{1}{7^2}+\frac{1}{8^2}-\ldots
\right)^{-1}\approx0.853\,
$$
in the case of packing, and
$$ D_3(\infty)=\left(
1-\frac{1}{2^2}+\frac{1}{4^2}-\frac{1}{5^2}+\frac{1}{7^2}-\frac{1}{8^2}+\ldots
\right)^{-1}\approx1.280
$$
in the case of covering (see also {\sc{Zeitler}} \cite{Zeitler} for the
case of covering). Of course, these densities are the same as the
density of the respective horoballs in the asymptotic tetrahedral cells of
the dual tiling $\{3,3,6\}$, that is, the bound proved by B\"or\"oczky for
packings and the still conjectured corresponding tetrahedral density bound
for coverings. Also, we have ${\lim\limits_{r\to\infty} d_3(r)=d_3(\infty)}$
and ${\lim\limits_{r\to\infty} D_3(r)=D_3(\infty)}$. {\sc Florian} (see
{\sc B\"or\"oczky} and {\sc Florian} \cite{BoroczkyFlorian}) showed that
in hyperbolic space the tetrahedral density bound $d_3=d_3(r)$ is a strictly
increasing function of $r$. Therefore, for an arbitrary packing of congruent
balls of finite or infinite radius in hyperbolic space, the density of
each ball in its Dirichlet cell is at most $d_3(\infty)$.  We can therefore
say that in 3-dimensional hyperbolic space the packing of the balls
inscribed in the cells of the tiling  $\{6,3,3\}$ is the densest among all
packings with congruent balls. It is conjectured that $d_n(r)$ is a
strictly increasing function of $r$ for all $n$. {\sc{Marshall}}
\cite{Marshall} gave a partial verification of this conjecture by proving
that it holds for sufficiently large values of $n$. {\sc{Kellerhals}}
\cite{Kellerhals} proved that $d_n(r)>d_{n+1}(r)$ for all $r>0$.

{\sc{Szirmai}} \cite{Szirmai05,Szirmai07,Szirmai18} and {\sc{Kozma}} and
{\sc{Szirmai}} \cite{KozmaSzirmai12,KozmaSzirmai15,KozmaSzirmai20}
determined the minimum density of certain classes of periodic packings
of horoballs in low dimensions.
\section{Hyperspheres}
Besides spheres of finite radius and horospheres, there is a third type of
``sphere" in hyperbolic space, namely hyperspheres. In $n$-dimensional
hyperbolic space a {\it{hypersphere}} is the set of points of the space
that are at the same distance $r$ from an $(n-1)$-dimensional hyperplane.
A hypersphere consists of two disjoint surfaces that bound a connected
component of the space called the {\it{hyperball of radius $r$}}. The study
of packing and covering by hyperballs was initiated by {\sc{Vermes}}.
In \cite{Vermes79} he investigated packings of congruent hypercircles
of radius $r$, and gave an upper bound for the local density of the
hypercircles in the cells of the dual subdivision of the plane into
Dirichlet cells. His bound is sharp for all values of
$r$, increases monotonously, and its limit at infinity is $3/\pi$, the
density of the densest packing of horocircles. Unaware of Vermes's result,
{\sc{Marshall}} and {\sc{Martin}} \cite{MarshallMartin03} discovered the
same bound.

{\sc{Przeworski}} \cite{Przeworski13} proved an upper bound for the density
of packings of congruent hyperballs. His bound is an analogue of the simplex
bound for ball packings. In the proof he uses {\it{Delone cells}} for the
base-hyperplanes of the hyperballs. These cells, which he introduced in
his paper \cite{Przeworski12}, are truncated ultraideal simplices. The
truncation faces are the base-hyperplanes of the hyperballs, and they are
perpendicular to all non-truncation faces that intersect them. The maximum
density of the hyperballs in the truncated simplices occurs for a regular
simplex. Hence he gets the following: In $H^n$, let $d_n(r)$ denote the density of
$n+1$ pairwise touching hyperballs in the truncated simplex bounded by the
base hyperplanes of the hyperballs and the hyperplanes orthogonal to the
base planes though the touching points of the hyperballs. Then the density of
any packing of hyperballs of radius $r$ within every Delone cell, as well as
in every Dirichlet cell, is at most $d_n(r)$. The same result was proved by
{\sc{Miyamoto}} \cite{Miyamoto} in another context. For 3-dimensional
space an alternative decomposition into truncated tetrahedra was suggested
by {\sc{Szirmai}} \cite{Szirmai19}.

{\sc{Vermes}} \cite{Vermes81} described for every $r>0$ a class of regular
packings of hypercircles of radius $r$, calculated the density of each
of them and determined the minimum density for all values of $r$. It turned
out that the minimum density $\Theta(r)$ is an increasing function of $r$
with limit value $\sqrt{12}/\pi$ as $r\to\infty$. In \cite{Vermes88} he
considered general coverings of hypercircles of radius $r$ subject to the
condition that there exists a number $0<\underline{r}<r$ such that the
hypercircles of radius $\underline{r}$ around the same base lines form
a packing. This condition guarantees that no point of the plane is covered
by infinitely many members of the covering. For the density of coverings of the
plane by hypercircles of radius $r$ satisfying the above property,
{\sc{Vermes}} \cite{Vermes88} gave a lower bound that is sharp for the
thinnest regular coverings. The infimum of the density of coverings by
hyperspheres is $\sqrt{12}/\pi$, attained by the thinnest covering of
horospheres.

A hyperball is the parallel body of an $(n-1)$-dimensional hyperplane. It is
natural to consider packings of parallel bodies of lower dimensional planes
as well. {\sc{Marshall}} and {\sc{Martin}} \cite{MarshallMartin05} and
{\sc{Przeworski}} \cite{Przeworski04,Przeworski06a} proved upper bounds
for the density of packings of {\it{tubes}}, that is of parallel bodies of
lines in 3-dimensional hyperbolic space.
\section{Solid arrangements}
As one of the possible substitutes for the notions of densest packing and thinnest
covering, {\sc{L.~Fejes T\'oth}} \cite{FTL68b} defined solidity of packings and
coverings. We say that a circle packing (covering) is {\it solid} if no finite subset
of the set of circles can be rearranged so that the resulting packing (covering)
is not congruent to the original one. Roughly speaking, if we remove any finite number
of arbitrarily chosen circles from a solid circle packing (covering), and we
wish to place them again so that we still obtain a packing (covering), then they have
to be returned to their original places.

It follows from a theorem of {\sc{Imre}} \cite{Imre} that the circles inscribed
in the faces of the tiling $\{k,3\}$ ($k=2,3,\ldots$) form a solid packing, and
the circles circumscribed about the faces of the tiling $\{k,3\}$ ($k=3,4, \ldots$)
form a solid covering. Moreover, {\sc{L.~Fejes T\'oth}} \cite{FTL68b} conjectured
that the circles inscribed in, and the circles circumscribed about the faces of
every three-valent Archimedian tiling form a solid packing and a solid covering,
respectively. On the other hand, he conjectured that the face-incircles and
the fave-circumcircles of a more than trihedral uniform tiling never form a solid
arrangement. In many cases this has been confirmed, see {\sc{L.~Fejes T\'oth}}
\cite{FTL68b}, {\sc{Heppes}} \cite{Heppes92}, {\sc{Heppes}} and {\sc{Kert\'esz}}
\cite{HeppesKertesz}, {\sc{Florian}} \cite{Florian99,Florian00,Florian01a,
Florian01b,Florian07}, {\sc{Florian}} and {\sc{Heppes}} \cite{FlorianHeppes00}
and {\sc{G.~Fejes T\'oth}} \cite{FTG74}. In particular, both conjectures were
confirmed for the face-incircles of all spherical and Euclidean uniform tilings.
As a particularly appealing example we mention the ``football'', {\it{i.e.}},
the tiling of the truncated icosahedron $\{5,6,6\}$, whose face-incircles
form a solid packing. In some cases, {\it{e.g.}}, for the incircles of the
tilings $(8,8,4)$ and $(4,4,n)$ for $n\ge6$, in addition to the solidity of the
packing it is shown that the arrangement has greatest possible density among all
packings with any collection of circles of the given radii (see {\sc{Florian}}
\cite{Florian01a}, {\sc{Florian}} and {\sc{Heppes}} \cite{FlorianHeppes99},
{\sc{Heppes}} \cite{Heppes00,Heppes03b} and {\sc{Heppes}} and {\sc{Kert\'esz}}
\cite{HeppesKertesz}).

We emphasize the following special case of the above conjecture: Augmenting
the incircles of a regular (spherical, Euclidean or hyperbolic) tiling by
circles inscribed in the holes results in a solid circle packing. This indeed
occurs in the case of the Euclidean tilings $\{6,3\}$ and $\{3,6\}$, since the
inscribed circles in $\{6,3\}$ form a solid packing by themselves,
and augmenting the incircles of $\{3,6\}$ results in the incircles
of $\{6,3\}$.  For the case of the third regular Euclidean tiling, namely
$\{4,4\}$, the solidity of the corresponding tiling was
proved by {\sc{Heppes}} \cite{Heppes92} using the idea of weighted density.

The solidity of a packing $\mathcal{P}$ follows if we can show that in every packing
of circles with the given radii the local density in the Delone triangles cannot
exceed the density of the circles of $\mathcal{P}$ in the Delone triangles. For the
face incircles of the tiling $(4,4,8)$ this is not true. Heppes observed that
the solidity of $\mathcal{P}$ follows if some weights can be assigned to the
circles so that the above statement holds for the weighted density. Using
appropriate weights he could prove, besides the incircles of the tiling $(4,4,8)$,
the solidity of several other packings. The work of Heppes inspired further research.
Weighted density of packings has also been studied for its own sake (see {\sc{H\'ars}}
\cite{Hars90,Hars92}).

A notion weaker than solidity was introduced and investigated by {\sc{ A. Bezdek,}}
{\sc{K. Bezdek}} and {\sc{Connelly}} \cite{BezdekABezdekKConnelly95,
BezdekABezdekKConnelly98}. A packing or a covering is {\it{uniformly stable}} if
there is an $\varepsilon>0$ such that no finite subset of the arrangement can be
rearranged so that each member is moved by a distance less than $\varepsilon$ and
the rearranged members, together with the rest, form a packing or covering,
respectively, different from the original arrangement. Using techniques from
rigidity theory they proved uniform stability of certain packings.

A solid circle packing is {\it{strongly solid}} if it remains solid even after
any one of the circles is removed. The
incircles of the faces of the tilings $\{4,3\}$ and $\{5,3\}$ are strongly solid.
L.~Fejes T\'oth conjectured that the incircles of the faces of the tiling
$\{p,3\}$ for $p\ge6$ form a strongly solid packing as well. {\sc{A.~Bezdek}}
\cite{BezdekA79} confirmed the conjecture for all $p\ge8$. The two remaining cases,
$p=6$ (in the Euclidean plane) and $p=7$ (in the hyperbolic plane) remain unsolved.
A partial result supporting the case $p=6$ of the conjecture is due to
{\sc{B\'ar\'any}} and {\sc{Dolbilin}} \cite{BaranyDolbilin}. They proved that the
packing obtained by removing one circle from the densest lattice packing of unit
circles is uniformly stable with $\varepsilon=1/40$. Another result supporting
the conjecture was given by {\sc{Heppes}} \cite{Heppes94}. He proved that the
hexagonal tiling $\{6,3\}$ is strongly translationally solid. This means that in the packing
of hexagons arising by omitting one tile from the hexagonal tiling, every rearrangement
of a finite number of hexagons by translations results in a packing congruent to the
original one.

{\sc{L.~Fejes T\'oth}} \cite{FTL80b} strengthened Bezdek's result in the
following sense. For any packing of congruent circles, a packing obtained
from it by removing a finite number $k$ of circles is called a {\it{$k$-truncation}}
of the packing. We say that a circle packing is {\it{solid of order $k$}}
if every $k$-truncation has the property that when finitely many
additional circles are removed from it and then placed back in any place where
there is room for them, the resulting packing is again some $k$-truncation of
the original packing. The {\it grade of saturation} of a packing of congruent
circles is the maximum number $g$ such that every circle congruent to, but not
identical with a circle of the packing intersects at least $g$ circles of the
packing.  {\sc{L.~Fejes T\'oth}} \cite{FTL80b} proved that for $p\ge8$, the order of solidity
of the face-incircles of $\{p,3\}$ equals the grade of saturation minus $2$.
In particular, it follows that the order of solidity of these packings becomes
arbitrarily large as $p\to\infty$.
\section{Completely saturated packings and completely reduced coverings}
{\sc{G.~Fejes T\'oth, G.~Kuperberg}} and {\sc{W.~Kuperberg}} \cite{FTGKuperbergKuperberg}
introduced the notion of a {\it $k$-saturated packing} $(k=1,2,\ldots)$ as a
packing with congruent copies of a set $K$, such that deleting $k-1$ members of
the packing never creates a void large enough to pack in it $k$ copies of $K$.
Similarly, a covering with congruent copies of $K$ is {\it $k$-reduced} if
deleting $k$ members of the covering always creates a void too large to be
covered by $k-1$ copies of $K$. A packing that is $k$-saturated for every $k$ is
{\it completely saturated}. {\it{Completely reduced}} coverings are defined similarly.
Obviously, in Euclidean spaces a completely saturated packing with congruent copies
of a convex body $K$ must have density $\delta(K)$; similarly, a completely reduced
covering with congruent copies of $K$ must be of density $\vartheta(K)$. This
suggested that the notions of completely saturated packings and completely reduced
coverings could serve in hyperbolic spaces as a substitute for packings and
coverings of maximum and minimum density, respectively. The question of existence of
completely saturated packings and completely reduced coverings for every convex
body turned out to be non-trivial by any means, and it was answered affirmatively
for Euclidean spaces in {\sc{G.~Fejes T\'oth, G.~Kuperberg}} and {\sc{W.~Kuperberg}}
\cite{FTGKuperbergKuperberg} and for hyperbolic spaces by {\sc{Bowen}} \cite{Bowen03a}.
The notion of completely reduced coverings found applications in approximation
theory, see {\sc{Hinrichs}} and {\sc{Richter}} \cite{HinrichsRichter04}.
\section{A probabilistic approach to optimal arrangements and their density}
{\sc{Bowen}} and {\sc{Radin}} \cite{BowenRadin03,BowenRadin04} proposed a
probabilistic approach to define optimal arrangements and their density,
especially useful in hyperbolic geometry. Their main idea
is that with a properly defined probability measure on the set of all packings
with copies of a body $K$, the density of a specific packing is the same as
the probability that its randomly chosen congruent copy contains the origin.
Below we sketch some details.

Let $\mathbb{S}$ denote an $n$-dimensional space of constant curvature, namely the
Euclidean $n$-space, the $n$-sphere, or the $n$-dimensional hyperbolic space.
Instead of studying individual packings in $\mathbb{S}$, Bowen and Radin consider
the space $\Sigma_K$ consisting of all packings of $\mathbb{S}$ by congruent
copies of a body $K$. A suitable topology derived from the Hausdorff metric on
$\Sigma_K$ is introduced which makes $\Sigma_K$ compact and makes the natural
action of the group ${\mathcal{G}}$ of rigid motions of $\mathbb{S}$ on $\Sigma_K$
continuous. We consider Borel probability measures on $\Sigma_K$ invariant under
${\mathcal{G}}$. For such an invariant measure $\mu$ the {\it density} of $\mu$,
$d(\mu)$, is defined as $d(\mu)=\mu(A)$, where $A$ is the set of packings
${\mathcal{P}}\in\Sigma_K$ for which the origin of $\mathbb{S}$ is
contained in some member of $\mathcal{P}$. It follows easily from the invariance
of $\mu$ that this definition is independent from the choice of the origin.

A measure $\mu$ is {\it ergodic} if it cannot be expressed as the positive
linear combination of two other invariant measures. The relationship between density
of measures and density of packings is established by the following theorem.

Suppose that $\mu$ is an ergodic invariant Borel probability measure on
$\Sigma_K$.  If a packing $P$ is chosen $\mu$-randomly, then with probability
$1$, for every $p\in {\mathbb{S}}$,
$$
d(\mu)=\lim_{\lambda\to\infty}\frac{1}{V(B_\lambda(p))}\sum_{P\in{\mathcal{P}}}V(P\cap{(B_\lambda(p))}.
$$

The {\it packing density} $\delta(K)$ of $K$ can then be defined as the
supremum of $d(\mu)$ for all ergodic invariant measures on $\Sigma_K$.  A
packing ${\mathcal{P}}\in\Sigma_K$ is {\it optimally dense} if the closure of its
orbit under ${\mathcal{G}}$ is the support of an ergodic invariant measure whose
density reaches this supremum.

It is shown in {\sc{Bowen}} and {\sc{Radin}} \cite{BowenRadin03,BowenRadin04} that,
for every $K$, an ergodic invariant measure $\mu$ with $d(\mu)=\delta(K)$ exists
whose support contains a set of full $\mu$-measure of optimally dense packings.
The existence of completely saturated packings with copies of $K$ was proved in the
same way as for optimally dense packings. In fact, with probability $1$, a
$\mu$-randomly chosen optimally dense packing is completely saturated.

{\sc{Bowen}} and {\sc{Radin}} \cite{BowenRadin03,BowenRadin04} proved several
statements justifying that they proposed a workable notion of optimal
density and optimally dense packings. In particular, in the Euclidean case,
that is ${\mathbb{S}} = {\mathbb{E}}^n$, the Bowen-Radin notion of $\delta(K)$
coincides with the corresponding traditional notion of the packing density of
$K$.

The probabilistic approach of Bowen-Radin can be naturally applied to
coverings, or more generally, to locally finite arrangements of congruent
copies of $K$.  The definition of density in such setting, however, requires a
modification: the ``measure'' $\mu(A)$ should be replaced by another quantity
that takes into account the multiplicity with which portions of the areas of
the union of $A$ are covered. More precisely, $\mu(A)$ should be replaced with a
combination of weighted measures, assigning the weight of $w$ to the regions
in $\mathbb{S}$ consisting of points covered exactly $w$ times by the given arrangement
of copies of $K$.

The advantage of this probabilistic approach is that it focuses on periodic
packings and neglects pathological packings such as B\"or\"oczky's example.
Another important feature of this approach is that we get bounds for the
packing density defined in this framework by considering the local density in
the Dirichlet or the Delone partitions.

Concerning packings of balls in hyperbolic $n$-space, it was shown by {\sc{Bowen}}
and {\sc{Radin}} \cite{BowenRadin03} that there are only countably many radii
admitting an optimally dense periodic packing of balls. Thus, for most radii $r$,
no periodic packing is densest. For the hyperbolic plane, {\sc{Bowen}}
\cite{Bowen03b} proved that the packing density of circles of radius $r$
is a continuous function of $r$, and it is the supremum of densities of periodic
packings.

\small{
\bibliography{pack}}